\documentclass{amsart}
\usepackage{amsmath,amssymb,amsfonts, amsthm}
\usepackage{filecontents}

\everymath{\displaystyle}

\newtheorem{theorem}{Theorem}
\newtheorem{definition}[theorem]{Definition}

\newtheorem{proposition}[theorem]{Proposition}
\newtheorem{corollary}[theorem]{Corollary}
\newtheorem{lemma}[theorem]{Lemma}

\def\R{\mathbb{ R}}

\def\N{\mathbb{N}}
\def\Psi{\varPsi}
\def\R{\mathbb{R}}
\def\A{{\mathcal P}_n}

\def\fb{\overline f}
\def\hb{\overline h}
\def\B{\mathcal B}

\begin{document}

\title[Sparse Hamburger Moment Sequences]{Sparse Hamburger Moment Sequences}



\author{Saroj Aryal}
\address[Saroj Aryal]
{Department of Mathematics\\
Montana State University\\
Billings, MT 59101}
\email[Saroj Aryal]{saroj.aryal@msubillings.edu}

\author{Hayoung Choi}
\address[Hayoung Choi and Farhad Jafari]
{Department of Mathematics\\
University of Wyoming\\
Laramie, WY 82071}
\email[Hayoung Choi]{hchoi2@uwyo.edu}

\author{Farhad Jafari*}
\email[Farhad Jafari]{fjafari@uwyo.edu}
\thanks{*: Corresponding author}
\thanks{Presented to the society at the memorial special session of the AMS regional meeting at the University of Memphis in October 2015.}

\subjclass[2010]{44A60, 13B30}

\date{}

\dedicatory{Dedicated to the memory of James Jamison}


\begin{abstract}
Putinar and Vasilescu \cite{putinar} have given an algebraic characterization of Hamburger moment sequences in several variables. In this paper we study some sparse moment subsequences of Hamburger moment sequences and consider the problem of completion of these moment subsequences. 
\end{abstract}
\maketitle

\section{Introduction.}
Let $\N_0$ be the set of nonnegative integers and $\N_0^n$ be the set of all multi-indices 
$\alpha =  (\alpha_1, \alpha_2, \ldots, \alpha_n)$,
$\alpha_j \in \N_0$ for all $j=1,2,\ldots,n$.  
A $n$-multisequence $(\gamma_\alpha)_{\alpha\in\N_0^n}$ is said to be \emph{positive semidefinite} if 
\begin{equation}\label{multipositive}
\sum_{\alpha,\beta\in S}\gamma_{\alpha+\beta} h_\alpha \hb_\beta \geq 0
\end{equation}
for any finite set $S\subset \N_0^n$ and complex-valued $n$-multisequence $(h_\alpha)_{\alpha \in S}$. We denote the set of positive semidefinite $n$-multisequences by {$\mathcal M_n$}. 

Let $ \A $ be the algebra of all complex-valued polynomials in $ n $-variables and  
$\B $ be the set of
nonnegative Borel measures $\sigma$ on $\R^n$ such that
\begin{equation}
\int_{\R^n}|x|^\alpha d\sigma(x)<\infty \quad \text{for all } \alpha \in\N_0^n.
\end{equation}
If $ L $ is a complex-valued linear map on $ \A $ such that $ L (1) > 0 $ and $ L(f\fb) \geq 0 $ ($ > 0 $) for every $ f \in \A$, then $ L $ is called positive semidefinite (definite). 
$ L $ is said to be a {\it moment map} if there exists a nonnegative Borel measure $\sigma \in \B$ such that $ L(f) = \int f d \sigma $ for all $ f \in \A $. Generalization of the Hamburger moment problem to higher dimensions is to find necessary and sufficient conditions for a real-valued multisequence $(\gamma_\alpha)_{\alpha\in\N_0^n}$ such that if $ L ( x^\alpha) = \gamma_\alpha $, then $ L $ is a moment map, i.e. there exists a nonnegative Borel measure $\sigma \in \B $ such that
\begin{equation}\label{multimomentproblem}
\gamma_\alpha=\int_{\R^n} x^{\alpha}d\sigma(x) \quad \textit{for all } \alpha\in \N_0^n,
\end{equation}
 If we denote the set of moment $n$-multisequences by $ M_n $, Hamburger's theorem states that $\mathcal M_1=M_1$, see \cite{akhiezer} for instance. For $ n > 1 $, it is easy to see $M_n \subset \mathcal M_n$. If the real-valued sequence $(\gamma_\alpha)_{\alpha\in\N_0^n}$ has the representation \eqref{multimomentproblem} for some nonnegative Borel measure $\sigma \in \B$, then 
\begin{align*}
\sum_{\alpha,\beta\in S}\gamma_{\alpha+\beta}h_\alpha h_\beta&=
\sum_{\alpha,\beta\in S} \int_{\R^n}x^{\alpha+\beta}\mathrm{d}\sigma(x) h_\alpha h_\beta\\
&=\int_{\R^n}\sum_{\alpha,\beta\in S} h_\alpha x^\alpha h_\beta x^\beta \mathrm{d}\sigma(x)\\
&=\int_{\R^n}\left[\sum_{\alpha\in S} h_\alpha x^\alpha\right]^2 \mathrm{d}\sigma(x)\geq 0,
\end{align*}
for all finite sets $S\subset \N_0^n$ and all $n$-multisequences $(h_{\alpha})_{\alpha\in S}$ in $\R^n$. Berg et al. \cite{berg} and Schm\"udgen \cite{schmudgen} independently proved that $\mathcal M_n\setminus M_n$ is nonempty. In a seminal paper, Putinar and Vasilescu \cite{putinar} provided algebraic conditions for a multisequence to be a Hamburger moment multisequence. In this paper, we consider the problem of identifying the subsequences of moment multisequences that are moment multisequences and consider the completion of these subequences. This work relies heavily on the techniques developed by Putinar and Vasilescu. We specialize those results by looking at the subideals of the algebras of rational functions containing $ \A $, and study the extension of the positive linear maps from these subalgebras to the full algebra.

In this work we identify a large class of subsequences of moment multisequences which are moment multisequences themselves. In Section 3 we show that while  the
moment subsequences of the form $ (\gamma_{d\alpha})_{\alpha \in \N_0^n} $ are completable, earlier results in one dimension show that moment subsequences may or may not
be completable in general. In Section 4 we show that there are subsequences of a moment sequence which are not moment sequences of the form \eqref{arith_sub}. 

It is useful to state the result from \cite{vasilescu} showing the additional required condition on $(\gamma_\alpha) \in\mathcal M_n$ to guarantee that it belongs to $M_n$. An equivalent result can be found in \cite{putinar}.
\begin{theorem}\label{vasilescu}
An $n$-multisequence $(\gamma_\alpha)_{\alpha \in \N_0^n} \in M_n  $ if and only if there exists a positive semidefinite 2n-multisequence $(\delta_{\alpha, \beta})_{\alpha, \beta\in \N_0^n } $ such that
\begin{align}
\label{delta1}\delta_{\alpha, 0} &= \gamma_\alpha \quad \text{and}\\
\label{delta2}\delta_{\alpha, \beta} &= \delta_{\alpha, \beta+ e_j} + \delta_{\alpha+ 2 e_j, \beta + e_j}, \quad 1 \leq j \leq n.
\end{align}
Furthermore, $\gamma$ is determinate if and only if the corresponding $\delta$ is unique.
\end{theorem}
We shall say that a multisequence $(\delta_{\alpha,\beta})\in\mathcal M_{2n}$ has the \textit{P-V property} if \eqref{delta1} and \eqref{delta2} hold. For the sake of notational convenience, some of the results will be stated just in two dimensions. Unless stated otherwise, these results can be trivially extended to higher dimensions.

\section{Main Result}
Given a moment multisequence $ (\gamma_\alpha)\in  \mathcal M_n$, define its subsequences
$(\tilde\gamma_\alpha)$ as
\begin{equation}\label{arith_sub}
\tilde\gamma_\alpha=\gamma_{d\alpha+\ell},
\end{equation}
where
$d=(d_1,d_2,\cdots,d_n)\in\N_0^n$, $\ell=(\ell_1,\ell_2,\cdots,\ell_n)\in2\N_0^n$, and $d\alpha=(d_1\alpha_1, \cdots,d_n\alpha_n).$
From the definition of positivity (\ref{multipositive}), it follows immediately that
\begin{proposition}\label{2dim_iff}
If $(\gamma_\alpha)\in\mathcal M_n$, then $(\tilde \gamma_\alpha)\in\mathcal M_n$. 
\end{proposition}

Now we are able to state our main result.
\begin{theorem}\label{main_multi}
If $(\gamma_\alpha)\in M_n$, then the subsequence $(\tilde \gamma_\alpha) \in M_n$.
\end{theorem}

To prove this theorem, we begin by proving a variation of Vasilescu's Theorem \cite{vasilescu}, viz,  
\begin{theorem}\label{var_vasilescu}
A $n$-multisequence $(\gamma_\alpha) \in M_n$ if and only if there exists $(\delta_{\alpha,\beta})\in\mathcal M_{2n}$ such that
\begin{align}
\label{delta111}\delta_{\alpha, 0} &= \gamma_\alpha\ \text{and}\\
\label{delta222}\delta_{\alpha, \beta} &= \delta_{\alpha, \beta+ e_j} + \delta_{\alpha+ 2 d_j e_j, \beta + e_j} \quad \text{for all }1 \leq j \leq n.
\end{align}
\end{theorem}
Here $d_j\in\N_0^n$ and $ (e_j) $ are the standard basis vectors in $ \R^n $. To prove this theorem, we will need a few lemmas. Let $\mathcal H$ be a complex Hilbert space whose scalar product will be denoted by $(*,*)$. If $ S $ is a linear map on $\mathcal H$, $S$ defines a linear subspace $D(S)\subseteq \mathcal H$ and has range $R(S)$. 

\begin{definition}
Let $T$ be a symmetric operator in $\mathcal H$. Then the Cayley transform of $T$ is an operator $U$ defined as
\begin{equation}
U=(T-i)(T+i)^{-1}.
\end{equation}
\end{definition}
\begin{lemma}\label{vasilesculemma2.0}
Let $U$ be the Cayley transform of a symmetric operator $T$ in $\mathcal H$. Then $U$ is unitary if and only if $T$ is self-adjoint.
\end{lemma}
\begin{proof}
See \cite[13.19]{rudin_FA}.
\end{proof}

\begin{lemma}\label{vasilesculemma2.1}
Let $S$ be a symmetric densely defined operator in $\mathcal H$. If the sets $R(S\pm i)$ are dense in $\mathcal H$, then the canonical closure of $S$ is a self-adjoint operator. \end{lemma}
\begin{proof}
Let $A$ be the canonical closure of $S$, which is also a symmetric operator. Since
\begin{equation}
\|(A\pm i)x\|^2=\|Ax\|^2+\|x\|^2,\ \ \ x\in D(A),
\end{equation}
then $R(A\pm i)$ are closed subspaces of $\mathcal H$. Since $R(A\pm i)\supseteq R(S\pm i)$ and $R(S\pm i)$ are dense in $\mathcal{H}$ by hypothesis, we have $R(A\pm i)=\mathcal H$. Let $V$ be the Cayley transform of $A$. Since $D(V)=R(A+i), R(V)=R(A-i)$, the operator $V$ is unitary. By Lemma \ref{vasilesculemma2.0}, $A$ must be self-adjoint.
\end{proof}
The following core lemma is a minor variation of Lemma 2.2 in \cite{vasilescu}. 
\begin{lemma} \label{vasilesculemma2.2}
Let $\theta_j(x)=\frac{1}{1+x_j^{2d_j}}$, $j=1,\ldots, n$, $x=(x_1,\ldots,x_n)\in\R^n$, $d=(d_1,\ldots,d_n)\in\N_0^n$. Let $\mathcal{R}_\theta$ be the algebra of rational functions generated by $\A $ and $\theta=(\theta_1,\ldots,\theta_n)$. 
Let $\rho:\mathcal{P}_{2n} \to \mathcal{R}_\theta$ defined as $\rho:p(x,s)\mapsto p(x,\theta(x))$. Then $\rho$ is a surjective unital algebra homomorphism, whose kernel is the ideal generated by the polynomials $\tau_j(x,s)=s_j\big(1+x_j^{2d_j}\big)-1$, $j=1,\ldots,n$, $s=(s_1,\cdots, s_n)\in\R^n$.
\end{lemma}
\begin{proof}
The fact that $\rho$ is a surjective unital algebra homomorphism is obvious. To determine the kernel of $\rho$, let $p\in \mathcal{P}_{2n}$ be a polynomial such that $p(x,\theta(x))=0$ for all $x\in\R^n.$ Then
\begin{equation}
p(x,s)=\sum_{\beta\in\N_0^n}p_\beta(x)s^\beta,
\end{equation}
with $p_\beta\in \A \backslash \{0\}$ only for a finite number of indices $\beta$. Then we have
\begin{align}\label{p_minus_p}
\nonumber p(x,s)&=p(x,s)-p(x,\theta(x))\\
\nonumber &=\sum_{\beta\neq 0}p_\beta(x)\left(s^\beta-\theta(x)^\beta\right)\\
&=\sum_{1\leq j\leq n}(s_j-\theta_j(x))q_j(x,s,\theta(x)),
\end{align}
where $q_j$ are polynomials. Let $a_j=\max_{1\leq j\leq n}\{\beta_j;p_\beta\neq0\}$, and define polynomials $t(x)$ as
\begin{equation}\label{poly_prod}
t(x)=\prod_{1\leq j\leq n}\left(1+x_j^{2d_j}\right)^{a_j}.
\end{equation}
Then from equations \eqref{p_minus_p} and \eqref{poly_prod}, we see that
\begin{equation}\label{tp}
t(x)p(x,s)=\sum_{1\leq j\leq n}(s_j(1+x_j^{2d_j})-1)\phi_j(x,s),
\end{equation}
where $\phi_j\in \mathcal{P}_{2n}$ for all $j=1,\cdots n$. If $a_j=0$ for all $j$, then $p(x,s)=p_0(x)=p(x,\theta(x))=0$. Therefore, without loss of generality, assume that $a_j\neq0$ for some $j$. Notice that the polynomials $t, \tau_j$ have no common zero in $\mathbb{C}^{2n}$. By a special case of Hilbert's Nullstellensatz, see \cite[Section 16.5]{waerden} for instance, there are polynomials $\tilde t, \tilde\tau_j\in \mathcal{P}_{2n}$ such that
\begin{equation}\label{ttildet}
t\tilde t+\sum_{1\leq j\leq n}\tau_j\tilde\tau_j=1.
\end{equation}
If we multiply \eqref{ttildet} by $p$, and use \eqref{tp}, we obtain
\begin{equation}\label{final_p}
p=\sum_{1\leq j\leq n}\tau_j(\phi_j\tilde t+\tilde\tau_jp).
\end{equation}
Thus kernel of $\rho$ is the ideal generated by the polynomials $\tau_j$.
\end{proof}
By the fundamental isomorphism theorem, $ \mathcal{P}_{2n}/{\ker} \rho $ is isomorphic to $ \mathcal{R}_\theta $. In particular, there is a one-to-one correspondence between these (commutative) algebras with units.
Now we are ready to prove Theorem \ref{var_vasilescu}.

\begin{proof}[Proof of Theorem \ref{var_vasilescu}]
Suppose $ (\gamma_\alpha)_{\alpha \in \N_0^n}  $ is a moment multisequence and hence has a representing measure $\sigma\in\B_n$. Define
\begin{equation}
\delta_{\alpha,\beta}=\int_{\R^n}x^{\alpha}\theta(x)^\beta d\sigma(x), \alpha, \beta \in \N_0^n,
\end{equation}
where
\begin{equation}\label{thetaj}
\theta_j(x)=\frac{1}{1+x_j^{2d_j}},~ 1\leq j\leq n, \text{ and } \theta=(\theta_1,\theta_2,\ldots,\theta_n).
\end{equation}
To see that $(\delta_{\alpha,\beta})$ is a positive sequence, consider a finite set $S\subseteq \N_0^n$ and a sequence of real numbers $(h_{\alpha,\beta})_{\alpha,\beta\in S}$. Then we have
\begin{align*}
\sum_{\alpha_1,\beta_1,\alpha_2,\beta_2\in S} \delta_{\alpha_1+\alpha_2, \beta_1+\beta_2}h_{\alpha_1,\beta_1}h_{\alpha_2,\beta_2}
&=\sum_{\alpha_1,\beta_1,\alpha_2,\beta_2\in S} \int_{\R^n}x^{\alpha_1+\alpha_2}\theta^{\beta_1+\beta_2}d\sigma(x)h_{\alpha_1,\beta_1}h_{\alpha_2,\beta_2}\\
&=\int_{\R^n}\sum_{\alpha_1,\beta_1\in S}\left[x^{\alpha_1}\theta^{\beta_1}h_{\alpha_1,\beta_1}\right]^2d\sigma(x),
\end{align*}
which is non-negative. Thus $(\delta_{\alpha,\beta})$ is positive.

Clearly, $(\delta_{\alpha,\beta})$ satisfies \eqref{delta111} as $$\delta_{\alpha,0}=\int_{\R^n}x^\alpha \theta(x)^0d\sigma(x)=\gamma_\alpha.$$ Similarly, for all $\alpha,\beta\in\N_0^n,~ j=1,\cdots,n$, we have
\begin{align*}
0 = &\int_{\R^n}\left(\theta_j(1+x_j^{2d_j})-1\right)x^\alpha \theta^\beta d\sigma(x)\\
= &\int_{\R^n}\left(\theta_j x^\alpha\theta^\beta+\theta_j x_j^{2d_j}x^\alpha\theta^\beta-x^\alpha\theta^\beta\right) d\sigma(x)\\
\implies &\delta_{\alpha,\beta+e_j}+\delta_{\alpha+2d_j e_j,\beta+e_j}=\delta_{\alpha,\beta},
\end{align*}
which establishes \eqref{delta222}.

Conversely, to prove the sufficiency direction, consider an $n$-sequence $\gamma$ and assume that there exists a positive $2n$-sequence $\delta$ satisfying \eqref{delta1} and \eqref{delta2}. Let $\theta_j$ and $\theta$ be as defined in \eqref{thetaj}. Let $R_\theta$ be an algebra generated by $\mathcal{P}_n$ and $\theta$. Define a positive semi-definite map $\Lambda$ on $R_\theta$ as
\begin{equation}
\Lambda(r) = L_\delta (p) , \  \ p \in \mathcal{P}_{2n} , \ \ r(x)= p(x,\theta(x))\in R_\theta,
\end{equation}
where $L_\delta$ is the linear map associated with $\delta$. Notice that, by Lemma \ref{vasilesculemma2.2}, the algebra $R_\theta$ is isomorphic to the quotient $\mathcal{P}_{2n}/I_{\tau}$, where $I_\tau$ is the ideal generated in $\mathcal{P}_{2n}$ by the polynomials $\tau_j(x,s)=s_j \big( 1+x_j^{2d_j} \big) -1$, $j=1,2,\ldots, n.$ Condition \eqref{delta222} implies $L_\delta|I_\tau=0$. Therefore, the map $L$ which can be identified with the map induced by $L_\delta$ on the quotient $\mathcal{P}_{2n}/I_{\tau}$ is well defined on $R_\theta$ and is positive semidefinite. Now we use the machinary that was developed by Gelfand and Naimark in their proof of the spectral theory of self-adjoint operators. Define a sesquilinear form on $R_\theta$ by the equation
\begin{equation}\label{slinear}
(r_1,r_2)_\theta=L(r_1\overline{r_2}), \ \ r_1, r_2\in R_\theta.
\end{equation}
Let $\mathcal N=\{r\in R_\theta:L(r\overline r)=0\}$. Then \eqref{slinear} induces a scalar product $(*,*)$ on the quotient $R_\theta/\mathcal N$. Let $\mathcal H$ be the completion of the quotient $R_\theta/\mathcal N$ with respect to this scalar product. Define in $\mathcal H$ the operators
\begin{equation}
T_j(r+\mathcal N)=x_jr+\mathcal N,~r\in R_\theta,~ j=1,\ldots,n,
\end{equation}
which are symmetric and densely defined, with $D(T_j)= R_\theta/\mathcal N$ for all $j$. The operators $T_j$ satisfy the conditions of Lemma \ref{vasilesculemma2.1} for every $j$. Indeed, if $r\in R_\theta$ is arbitrary, then the functions $u_{\pm}(x)=(x_j\mp i\theta_j(x)r(x))$ are solutions in $R_\theta$ of the equations $(x_j\pm i)u_{\mp}(x)=r(x).$ This implies that the equalities $R(T_j\pm i)=D(T_j)$, and therefore using Lemma \ref{vasilesculemma2.1}, $T_j$ is essentially self-adjoint.

Let $A_j$ be the canonical closure of $T_j$. Now, we want to show that the operators $(i-A_1)^{-1},\cdots,(i-A_n)^{-1}$ mutually commute. Notice that the map $(i-T_j)^{-1}$ is well defined on $D=D(T_j)$ and invariant $j$, by the above argument. Moreover, the maps $(i-T_1)^{-1},\cdots,(i-T_n)^{-1}$ mutuallly commute on $D$. Since $A_j$ extends $T_j$, we have
$$
(i-A_j)((i-A_j)^{-1}-(i-T_j)^{-1})\xi=0, \ \ \xi\in D,
$$
implying
$$(i-A_j)^{-1}|D=(i-T_j)^{-1}.$$ 
Therefore, for all $j,k$ with $j\neq k$, we have
\begin{align*}
(i-A_j)^{-1}(i-A_k)^{-1}\xi
&=(i-T_j)^{-1}(i-T_k)^{-1}\xi\\
&=(i-T_k)^{-1}(i-T_j)^{-1}\xi\\
&=(i-A_k)^{-1}(i-A_j)^{-1}\xi,\ \ \xi\in D.
\end{align*}
Since $(i-A_1)^{-1},\ldots,(i-A_n)^{-1}$ are bounded and $D$ is dense, they mutually commute. Therefore, the self-adjoint operators $A_1,\ldots,A_n$ have a joint spectral operator measure $E$. Then the measure
\begin{equation}\label{sigmaE}
\sigma(*)=(E(*)(1+\mathcal N),1+\mathcal N),
\end{equation}
is a representing measure for the functional $\Lambda$. In what follows, we will show $\sigma$ is a representing measure for the given sequence $\gamma$.

Let $r(T)$ be the linear map on $D$ given by
$$
r(T)(f+\mathcal N)=rf+\mathcal N, \ \ r,f\in R_\theta.
$$
then we have $\theta(A)^\beta\supset\theta(T)^\beta$, where $\theta(A)^\beta$ is given by the functional calculus of $A$. This follows from the relations $\theta(A)^{-\beta}\supset\theta(T)^{-\beta}$ and $\theta(A)^{-\beta}(\theta(A)^{\beta}-\theta(T)^{\beta})=0$. 
Therefore,
\begin{align*}
\delta_{\alpha,\beta}&=(x^\alpha\theta(x)^\beta1,1)_\theta\\
&=(t^\alpha\theta(T)^\beta (1+\mathcal N),1+\mathcal N)\\
&=(A^\alpha\theta(A)^\beta(1+\mathcal N),1+\mathcal N)\\
&=\int_{\R^n}x^\alpha \theta^\beta d(E(x)(1+\mathcal N), 1+\mathcal N).
\end{align*}
Hence, using the condition \eqref{delta1} and equation \eqref{sigmaE}, we have
$$
\gamma_\alpha=\delta_{\alpha,0}=\int_{\R^n}x^\alpha \theta^0 d(E(x)(1+\mathcal N), 1+\mathcal N)=\int_{\R^n}x^\alpha d\sigma(x).
$$
\end{proof}

\begin{proof}[Proof of Theorem \ref{main_multi}]
Suppose $\tilde\gamma_\alpha=\gamma_{d\alpha+\ell}$ where $d=(d_1,\cdots,d_j)\in\N_0^n$ and $\ell\in2\N_0^n$. Since $(\gamma_\alpha) \in M_n$, there is a multisequence $(\delta_{\alpha,\beta})\in\mathcal M_{2n}$ satisfying the equations \eqref{delta111} and \eqref{delta222} in Theorem \ref{var_vasilescu}. Define another multisequence $(\tilde\delta_{\alpha,\beta})$ as
\begin{equation}
\tilde\delta_{\alpha,\beta}=\delta_{d\alpha+\ell,\beta}.
\end{equation}
It can be easily verified that $(\tilde\delta _{\alpha,\beta})\in\mathcal M_{2n}$ and satisfies \eqref{delta1}. Applying the condition of equation \eqref{delta222} to $\delta_{d\alpha+\ell,\beta}$ we have for every $j=1,\cdots,n$,
\begin{equation*}
\delta_{d\alpha+\ell,\beta}=\delta_{d\alpha+\ell,\beta+e_j}+\delta_{d\alpha+\ell+2de_j,\beta+e_j},
\end{equation*}
which can be re-written as
\begin{equation*}
\delta_{d\alpha+\ell,\beta}=\delta_{d\alpha+\ell,\beta+e_j}+\delta_{d(\alpha+2e_j)+\ell,\beta+e_j},
\end{equation*}
implying
\begin{equation}
\tilde\delta_{\alpha,\beta}=\tilde\delta_{\alpha,\beta+e_j}+\tilde\delta_{\alpha+2e_j,\beta+e_j},
\end{equation}
which is just equation \eqref{delta2} for $(\tilde\delta _{\alpha,\beta})$. Therefore, by Theorem \ref{vasilescu}, we conclude that $(\tilde\gamma_\alpha )\in M_n$.
\end{proof}

As a concluding remark to this section, we provide below another sufficient condition for the subsequence $(\tilde\gamma_\alpha)$ to be a moment sequence.

\begin{theorem}\label{suff_multi}
A subsequence $(\tilde\gamma_\alpha)$ as defined in \eqref{arith_sub} of a moment multisequence $(\gamma_\alpha)_{\alpha\in\N_0^n}$ having a representing measure $\sigma$ is a moment multisequence if
\begin{equation}\label{condition_multi}
\int_{\R^n}x^{d\alpha+\ell}\theta^\beta\frac{x_j^{2d_j}-x_j^{2}}{x_j^2+1}d\sigma(x)=0, \text{ for every } j=1,\cdots,n,~and~ \beta\in\N_0^n,
\end{equation}
where $\theta=(\theta_1,\cdots,\theta_n)$ and $\theta_j=\frac{1}{1+x_j^2}$.
\end{theorem}
\begin{proof}
 By Theorem \ref{vasilescu}, there exists a positive semidefinite  $ 2n $-multisequence $ (\delta_{\alpha, \beta}) $ with the P-V property that extends $ (\gamma_\alpha) $, where $(\delta_{\alpha,\beta})$ is defined as
 \begin{equation}
 \delta_{\alpha,\beta}=\int_{\R^n}x^\alpha\theta^\beta d\sigma(x).
 \end{equation}
 Now, define a $2n$-multisequence $(\tilde\delta_{\alpha,\beta})$ as
 \begin{equation}
 \tilde\delta_{\alpha,\beta}=\delta_{d\alpha+\ell,\beta}, \ \ \alpha,\beta\in\N_0^n.
 \end{equation}
 Clearly, $(\tilde\delta_{\alpha,\beta})$ is positive by Theorem \ref{2dim_iff}. It is sufficient to show that $(\tilde\delta_{\alpha,\beta})$ satisfies the P-V property, namely conditions \eqref{delta1} and \eqref{delta2} hold.
 Since
 \begin{equation}
 \tilde\gamma_\alpha=\gamma_{d\alpha}=\delta_{d\alpha, 0}=\tilde\delta_{\alpha,0},
 \end{equation}
 then condition \eqref{delta1} is satisfied. To prove condition \eqref{delta2}, we have
 \begin{align*}
\tilde\delta_{\alpha,\beta}
&=\delta_{d\alpha+\ell,\beta}\\
&=\delta_{d\alpha+\ell, \beta+ e_j} + \delta_{d\alpha+ \ell+2 e_j, \beta + e_j}\\
&=\int_{\R^n}\left(x^{d\alpha+\ell}\theta^{\beta+ e_j}+x^{d\alpha+ \ell+2 e_j}\theta^{\beta+ e_j}+x^{d\alpha+\ell}\theta^\beta\frac{x_j^{2d_j}-x_j^{2}}{x_j^2+1}\right)d\sigma(x) \\
&=\int_{\R^n}\left(x^{d\alpha+\ell}\theta^{\beta+ e_j}+x^{d(\alpha+2e_j)+ \ell}\theta^{\beta+ e_j}\right)d\sigma(x)\\
&=\delta_{d\alpha + \ell, \beta+ e_j} + \delta_{d\alpha+ 2d e_j+\ell, \beta + e_j}\\
&=\tilde\delta_{\alpha, \beta+ e_j} + \tilde\delta_{\alpha+ 2 e_j, \beta + e_j} \quad \text{for all }j\in \N_0.
 \end{align*}
Therefore, $(\tilde\delta) $ satisfies the P-V property and $ (\tilde\gamma) \in M_n$.
\end{proof}

\begin{corollary}\label{multi_shift}
If $(\gamma_\alpha)_{\alpha\in\N_0^n}\in M_n$, then $(\tilde\gamma_\alpha)$ defined as $\tilde\gamma_\alpha=\gamma_{\alpha+\ell}$ is in $M_n$.
\end{corollary}
\begin{proof}
Notice that the condition \eqref{condition_multi} in Theorem \ref{suff_multi} is satisfied when $d_j=1$ for all $j\in\N_0$.
\end{proof}
The subsequence identified in Corollary \ref{multi_shift} is particularly useful in the case when a block of data is missing. 

\section{Completions}
Given a positive map $ \tilde{ \Lambda} $ on the subalgebra 
$ \tilde{\mathcal R}_\theta $ determined by the ideal generated by the polynomials $ \tau_j  (x,s) = s_j (1+x_j^{2d_j}) -1 $, $ j = 1, \cdots, n $, one may ask if this map $ \tilde{\Lambda}$ extends to a positive map $ \Lambda $ on the larger algebra, $ {\mathcal R}_\theta $ determined by the ideal generated by $ \eta_j (x,s) = s_j (1+x_j^{2})-1 $, $ j=1, \cdots, n $. If so, by Corollary 2.6 in \cite{putinar}, there is a uniquely determined representing measure $ \tilde{\mu} $ in $ \R^n $ for $ \tilde{ \Lambda} $ such that the algebra 
$ \tilde{\mathcal R}_\theta $ is dense in $ L^2 (\tilde\mu) $. Since, without loss of generality, we may assume $ \tilde{\Lambda} (e) = 1 $ for the unit $ e$  in  $ \tilde{\mathcal R}_\theta $, by the Hahn-Banach theorem for positive functionals, the positive functional $ \tilde{\Lambda} $ extends to a positive functional $ \Lambda $ on $ {\mathcal R}_\theta $. Hence, again applying Corollary 2.6 in \cite{putinar}, there exists a $ \mu $ on $ \R^n $ such that $ \Lambda (f) = \int f d\mu $ for every $ f \in {\mathcal R}_\theta $, $ {\mathcal R}_\theta $ is dense in $ L^2 (\mu) $, and $ \Lambda $ agrees with $\tilde{\Lambda} $ on the smaller algebra. In particular, the partial sequence is completable. In general, since the extension of $ \tilde{\Lambda} $ to $ \Lambda $ is not unique, this extension will not be unique. Thus, we have 
\begin{theorem}\label{completion}
Let $(\gamma_\alpha)_{\alpha\in \N_0^n}$ be a partial sequence 
defined by
\begin{equation*}
\gamma_\alpha=
\left\{
\begin{array}{rl}
\textit{specified}& \text{if } \alpha \in d\N_0^n,\\
\textit{missing} & \text{if } \alpha \notin d\N_0^n
\end{array}\right.
\quad\textit{for all } \alpha \in \N_0^n.
\end{equation*}
If its subsequence 
$ ({\tilde{\gamma}_\alpha})_{\alpha\in \N_0^n}$  defined by
$$ \tilde{\gamma}_\alpha = \gamma_{d\alpha} \quad \text{for all } \alpha \in \N_0^n$$
is in $M_n$, then the partial sequence $(\gamma_\alpha)_{\alpha\in \N_0^n}$ 
has a completion with fully specified entries for all $\alpha \in \N_0^n$. 
\end{theorem}


\section{Some Counterexamples}

The above shows that if $(\gamma_\alpha)\in\mathcal M_n$ is arbitrary, then its subsequences $(\tilde \gamma_\alpha)\in\mathcal M_n$. However, for specific sequences $ (\gamma_\alpha) $, or even for those with various growth properties, the below examples show that no reasonable characterization of moment subsequences is possible.

\subsection*{Example 1}
We state below one of the main results from \cite{craven} without proof.

\begin{theorem}\label{strictpositivity}
Let $(s_k)_{k\in\N_0}$ be a sequence of positive numbers satisfying
\begin{equation}\label{craveneq}s_{k+1}s_{k-1}\geq c_0s_k^2,\end{equation}
where $c_0$ is the unique real root of $x^3-5x^2+4x-1$ approximately equal to $4.0796$. Then, for each positive integer $n$, the Hankel matrices $H=\{s_{i+j-2}\}_{1 \leq i,j \leq n}$ and $H_1=\{s_{i+j-1}\}_{1 \leq i,j \leq n}$ are strictly totally positive. Moreover, there is a nondecreasing function $\sigma$ with infinitely points of increase such that
\begin{equation}\label{stieltjesmoments}
s_n=\int_0^{\infty}x^nd\sigma(x) \quad \text{for all } n\in\N_0.
\end{equation}
\end{theorem}

Note that the integral representation \eqref{stieltjesmoments} says that $(s_k)$ is a Stieltjes moment sequence. However, the positivity of the Hankel matrix $H=\{s_{i+j-2}\}_{1\leq i,j\leq n}$ makes it a Hamburger moment sequence as well. Indeed, any nondecreasing measure in $[0,\infty)$ can be trivially extended to a nondecreasing measure in $\R$.

Consider a subsequence $(s_k)$ that satisfies the condition \eqref{craveneq} and hence is a Hamburger moment sequence using Theorem \ref{strictpositivity}. Let $(t_k)$ be its subsequence such that $t_0=s_0$, $t_1=s_1$, and $t_k=s_{k+1}$ for all $k\geq 2$. Notice that $(t_k)$ is not in an arithmetic pattern of $(s_k)$. But since
$$t_{k+1}=s_{k+2}\geq c_0\cfrac{s_{k+1}^2}{s_{k}}=c_0\cfrac{t_{k}^2}{t_{k-1}},$$
the sequence $(t_k)$ is a Hamburger moment sequence using Theorem \ref{strictpositivity}.

\subsection*{Example 2}
The Hilbert sequence $(s_k)_{k\in \N_0}$ given as $s_k=\dfrac{1}{k+1}$ is a positive sequence since determinant of its Hankel matrix is  
$$D_n=\left|
\begin{smallmatrix} \dfrac{1}{1} & \dfrac{1}{2} & \dfrac{1}{3} & \cdots & \dfrac{1}{n+1} \\\dfrac{1}{2} & \dfrac{1}{3} & \dfrac{1}{4} & \cdots & \dfrac{1}{n+2} \\\dfrac{1}{3} & \dfrac{1}{4} & \dfrac{1}{5} & \cdots & \vdots \\\vdots & \vdots & \vdots & \ddots & \vdots \\\dfrac{1}{n+1} & \dfrac{1}{n+2} & \cdots & \cdots & \dfrac{1}{2n+1}
\end{smallmatrix}
\right|=\prod_{i=0}^n\left[(2i+1)\binom{2i}{i}^2\right]>0$$
for all $n\in\N_0$.
So $(s_k)$ has a moment solution by Hamburger's theorem. Clearly the sequence $t_k=\left(\dfrac{1}{k+1}\right)^2$ is a subsequence of $(s_k)$. Then $(t_k)$ is a Schur square of $(s_k)$ and hence is a moment sequence (see \cite{horn}). But $(t_k)$ is not of the form \eqref{arith_sub}.

Of course, there is a deeper question in the background here. 
If $ P \subset \N_0^n $ is a  {\it pattern}, in Section 2 it is shown that if $ (\gamma_\alpha) $ is a moment multisequence, and $ P $ is a multidimensional arithmetic progression, then $ (\gamma_\alpha)_{\alpha \in P} $ is a moment multisequence. The converse of this theorem, would ask if $ (\gamma_\alpha) $ is any (arbitrary) moment sequence and $ (\gamma_\alpha)_{\alpha \in P} $ is a moment multisequence, is $ P $ necessarily an arithmetic progression?  The fact that this is false even for $ n=1 $ was shown in \cite{choijafari}.

\section*{Acknowledgements}
The authors wish to express their gratitude to the anonymous referee(s) for their
careful reading of this manuscript and for their valuable suggestions.  The suggested changes
have improved the paper significantly.  Also, we would like to thank the editors of this
volume for organizing a memorial conference in honor of James Jamison, and for dedicating this
issue to him.

\bibliographystyle{amsplain}
\bibliography{refs}

\end{document}